\documentclass[%
 reprint,
 amsmath,amssymb,
 aps,
]{revtex4-1}

\usepackage{subfigure}
\usepackage{graphicx}% Include figure files
\usepackage{dcolumn}% Align table columns on decimal point
\usepackage{bm}% bold math
\usepackage[colorlinks,linkcolor=blue,anchorcolor=blue,citecolor=blue]{hyperref}

\begin{document}

\preprint{APS/123-QED}
\title{Fractional compound Poisson processes with multiple internal states }

\thanks{This work was supported by the National Natural Science Foundation of China under Grant No. 11671182.}%

\author{Pengbo Xu}
 \email{xupb09@lzu.edu.cn}%Lines break automatically or can be forced with \\
\author{Weihua Deng}%
 \email{dengwh@lzu.edu.cn}
\affiliation{%
 School of Mathematics and Statistics, Gansu Key Laboratory of Applied Mathematics and Complex Systems, Lanzhou University, Lanzhou 730000, P.R. China
}%


\date{\today}% It is always \today, today,

\begin{abstract}
For the particles undergoing the anomalous diffusion with different waiting time distributions for different internal states, we derive the Fokker-Planck and Feymann-Kac equations, respectively, describing positions of the particles and functional distributions of the trajectories of particles; in particular, the equations governing the functional distribution of internal states are also obtained. The dynamics of the stochastic processes are analyzed and the applications, calculating the distribution of the first passage time and the distribution of the fraction of the  occupation time, of the equations are given.

% The applications of the equations are given, including
%\begin{description}
%\end{description}
\end{abstract}
\maketitle

\emph{Introduction.}--- Poisson process is the most fundamental stochastic process of renewal theory. The application of Poisson process naturally coming to our minds is in queueing theory to model the random events: the arrival of customers at a store or phone calls at an exchange \cite{Kleinrock:76}. Renewal process generalizes Poisson process for arbitrary holding times, being still independent identically distributed (i.i.d.) \cite{godreche:00}. If the probability density function (PDF) of the holding/waiting times between two subsequent events has the asymptotic behavior $\phi(\tau)\sim 1/\tau^{\alpha+1}, 0<\alpha<1$ when time is long enough \cite{laskin:00}, it is called fractional Poisson process.
The residence time statistics for $N$ fractional Poisson  processes are discussed in \cite{burov:00}, which reveal sharp transitions for a critical number of degrees of freedom $N$ and give application of detecting nonergodic kinetics from the measurements of many blinking chromophores.

We further generalize the renewal processes to have multiple internal states, where the holding times for different internal states are drawn from different distributions. The case of two internal states is considered in \cite{lowen:00,Godec:17} with applications, including trapping
in amorphous semiconductors, electronic burst noise, movement in systems with fractal boundaries,
the digital generation of $1/f$ noise, and ionic currents in cell membranes. The fractional Poisson processes with multiple internal states have a lot of potential applications, e.g., the particles moving in multiphase viscous liquid composed of materials with different chemical properties; Niemann, Barkai, and Kantz \cite{niemann:00} detailedly investigate a stochastic signal with multiple states, in which each state has an associated joint distribution for the signal's intensity and its holding time. This letter adopts the notations of \cite{niemann:00}.
If introducing jumps to the fractional Poisson processes, then we reach the fractional compound Poisson (fcP) processes \cite{meerschaert:00}. The jumps follow the fractional Poisson processes to arrive and, most of the time, the size of the jumps is random, with specified PDF. More specifically, let $\mathcal{N}(t)$ be the fractional Poisson process and $\{\xi_i, i=1,2,\ldots\}$ a sequence of i.i.d. random variables (jump lengths). Then $X(t)=\sum\limits_{i=1}^{\mathcal{N}(t)} \xi_i$ is the fcP process. The continuous time random walk (CTRW) with i.i.d.  power law waiting time \cite{cartea:00, metzler:00,Metzler:16,klafter:00} is a specific fcP process, being widely used to model various anomalous diffusions, e.g., mRNA molecules in living cells \cite{golding:00}, price fluctuation in financial market \cite{scalas:00}.

% and, most of the time, the size of the jumps is
%It's well known that the Poisson process is a counting process. And it's one of the most studied and applied processes. If we introduce a concept of the waiting time to the fractional Poisson processes. Then the waiting time has the asymptotic behavior $\phi(\tau)\sim 1/\tau^{\alpha+1}, 0<\alpha<1$ when time is long enough \cite{laskin:00}. Therefore, from this point of view Poisson process is the fundamental theory of the renewal process. Besides, if a process, denoted as ${X(t),t\geq0}$, can be represented as $X(t)=\sum_{i=1}^{N(t)}\xi_i$ where $\{N(t)\}$ is the fractional Poisson process, and $\{\xi_i, i=1,2,\ldots\}$ is a sequence of i.i.d. variables, then we call this process $X(t)$ is fractional compound Poisson(FCP) process. The FCP process with 2 alternating internal states is discussed in ref. \cite{lowen:00} and has many applications in quantum jumps and fluorescence signal from a quantum dot\cite{burov:00,godreche:00}.
%In this letter we mainly study the FCP process with several internal states to choose one of the distributions of the waiting time to obey.

In this letter, we focus on the fcP processes $X(t)$  with multiple internal states, i.e., $\mathcal{N}(t)$ of $X(t)$ is a fractional Poisson process with multiple internal states. Each internal state has an associated distribution of waiting time, but the distributions of jump lengths are all simply taken as normal distribution. We derive the Fokker-Planck equations \cite{metzler:00}, governing the PDF of positions of $X(t)$, the Feynman-Kac equations \cite{turgeman:00, carmi:00, cairoli:00}, describing the distribution of the functional \cite{satya:00} of the paths of $X(t)$, and the equations, characterizing the functional distribution of the internal states. From the Fokker-Planck equations, we obtain the evolution of the mean square displacement (MSD) for the process, showing that it strongly depends on the properties of the internal transition matrix \cite{feller:00} and sometimes the distribution of the initial position of the particles. The applications of the Feynman-Kac equations and the equations governing the distribution of the functional are given to calculate the distribution of first passage time \cite{carmi:00,redner:00,Deng:17} and the distribution of the fraction of the occupation time \cite{godreche:00}. Some properties of the distributions are obtained.

%we mainly obtain the Fokker-Planck equations and corresponding Feynman-Kac equations. From these equations we can obtain some important properties, such as mean square displacement(MSD) of this process, the distribution of the first passage time, and the distribution of the friction of occupation time and the total time. The final results rely on the transition matrix and the waiting time distributions. In this letter we take the distributions of jump length are the same normal distribution, thus we don't mention the affects of the jump length distributions. As for the relationship between the results and the initial distribution of the internal states, it depends on the transition matrix. More specific, if the transition matrix is irreducible, then the initial distribution has no concern with the final conclusions. On the other hand, when the transition matrix is reducible, the initial distribution will influence the results.

\emph{Model.}--- We consider the fcP processes with finite internal states, denoting their number as $N$.
%The internal states has many abstract meanings, such as 'on' and 'off'. The number of internal states is denoted as $N$.
%We will use the internal states to choose corresponding distributions of waiting time and jump length to obey. The specific way of choosing will be illustrated in the following.
The internal states determine the distributions of waiting times and
%For the convenient to state clearly, we would introduce some symbols first.
the transition of the internal states is described by a Markov chain with its transition matrix $M$; the dimension of $M$ is $N\times N$. The element $m_{ij}$ of the matrix $M$ represents the transition probability from state $i$ to state $j$.
%the possibility of the internal state transiting from $i$th to $j$th.
The bras $\big<\cdots\big|$ and kets $\big|\cdots\big>$ denote the row and column vectors,  respectively. From Ref. \cite{feller:00} one can see that for the ergodic or periodic chain there always exists an equilibrium distribution, denoted by $\big<\rm eq_M\big|$, with the property of $\big< {\rm eq_M}\big|M=\big<{\rm eq_M}\big|$. For the transition matrix $M$, as is well known, its largest eigenvalue is $1$ and the right eigenvector for this eigenvalue is the vector of 1s \cite{Walker:11}, denoted as $\big|\Sigma\big>$. That is $M\big|\Sigma\big>=\big|\Sigma\big>$. The second largest eigenvalue is strictly less than $1$ if $M$ is irreducible. It also obviously holds that $M^{T}\big|{\rm eq_M}\big>=\big|{\rm eq_M}\big>$ and $\big<\Sigma\big|M^{T}=\big<\Sigma\big|$.
%we can also find one of the right eigenvectors denoted by $\big|\Sigma\big>$ which every factor of this vector is 1. Thus we have $M\big|\Sigma\big>=\big|\Sigma\big>$.
We use the notation $\big<{\rm init}\big|$ to represent the initial distribution of the internal states.
%Apparently, we have $M^{T}\big|{\rm eq_M}\big>=\big|{\rm eq_M}\big>$ and $\big<\Sigma\big|M^{T}=\big<\Sigma\big|$.
Based on the CTRW, we define the waiting time distribution matrix $\Phi(t)= {\rm diag} (\phi^{(1)}(t),\phi^{(2)}(t),\ldots,\phi^{(N)}(t))$ and the jump length one $\Lambda(x)={\rm diag}(\lambda^{(1)}(x),\lambda^{(2)}(x),\ldots,\lambda^{(N)}(x))$.

The initial state of the stochastic process $X(t)$ is drawn from the initial distribution. After confirming the initial state, e.g., state $i$, the waiting time and jump length will be obtained from the distributions $\phi^{(i)}(t)$ and $\lambda^{(i)}(x)$, respectively. The new internal state is drawn from the distribution $M\big|{\rm init}\big>$. Then repeat the procedure.
%The particles choose which state to stay according to the initial distribution when the stochastic process start to proceed. After confirming the internal state, such as $i$th state, then we choose the waiting time and the jump length choose corresponding distributions $\phi^{(i)}(t)$ and $\lambda^{(i)}(x)$ to follow. And for the second step, the particles choose internal state according to the distribution $M\big|{\rm init}\big>$. Then repeat the procedures.
For the Markov chain, there are plenty of practically or theoretically important transition matrices \cite{feller:00}. The simple and representative one should be $\begin{pmatrix}
0 & 1\\
1 & 0\\
\end{pmatrix},$
signifying $2$ alternating internal states. The model will be, respectively, discussed for the case that $M$ is irreducible and the case that $M$ is reducible, i.e., the digraph of $M$ is strongly connected or not.

%its associated digraph is not strongly connected

\emph{Fokker-Planck equations.}---We use the notation $g^{(i)}(x,t)$, $i=1,2,\ldots,N$ to represent the PDF of finding the particle, at time $t$, position $x$ and internal state $i$. Let   $\big|G(x,t)\big>$ be the column vector $\{g^{(i)}(x,t),\,i=1,\ldots,N\}$. The initial condition for $\big|G(x,t)\big>$ is taken as $\delta(x)\delta(t)\big|{\rm init}\big>$. Similarly to the derivation of fractional Fokker-Planck equation from CTRW model with i.i.d. waiting times \cite{klafter:00,metzler:00}, we can obtain the equation in the Fourier-Laplace space
\begin{equation}\label{eq0}
\big|G(k,s)\big>=\frac{1}{s}\big[I-\Phi(s)\big]\big[I-M^{T}\Phi(s)\Lambda(k)\big]^{-1}\big|{\rm init}\big>.
\end{equation}
In this letter, we take the waiting time distributions as asymptotical power laws, i.e., in the Laplace space $\Phi(s)=I-\Psi(s)$, where $\Psi(s)={\rm diag}(B_{\alpha_1}s^{\alpha_1},\ldots,B_{\alpha_N}s^{\alpha_N})$, $0<\alpha_1,\ldots,\alpha_N<1$. As for the jump length distributions we choose all of them as Gaussian distribution, so in the Fourier space $\Lambda(k)=(1-\sigma^2k^2)I$. Thus from the above equations, %we have
%\begin{equation}\label{eq1}
%\begin{split}
%sM^T\big|G(k,s)\big>-\big|{\rm init}\big>=&-s\big(I-M^T+\sigma^2k^2M^T\big)\\
%&\Psi^{-1}(s)\big|G(k,s)\big>.
%\end{split}
%\end{equation}
taking the inverse Laplace and Fourier transformations leads to the Fokker-Planck equation with $N$ internal states
\begin{widetext}
\begin{equation}\label{eq2}
\begin{split}
M^T\frac{\partial}{\partial t}\big|G(x,t)\big>=&\big(M^T-I\big){\rm diag}\big(B_{\alpha_1}^{-1},\ldots,B_{\alpha_N}^{-1}\big)D^{{\rm diag}(1-\alpha_1,\ldots,1-\alpha_N)}_t\big|G(x,t)\big>\\
&+M^T{\rm diag}\big(K_{\alpha_1},\ldots,K_{\alpha_N}\big)D^{{\rm diag}(1-\alpha_1,\ldots,1-\alpha_N)}_t\frac{\partial^2}{\partial x^2}\big|G(x,t)\big>,
\end{split}
\end{equation}
\end{widetext}
where %$D^{{\rm diag}(1-\alpha_1,\ldots,1-\alpha_N)}_t$ is defined as
$
D^{{\rm diag}(1-\alpha_1,\ldots,1-\alpha_N)}_t={\rm diag}\Big(D_t^{1-\alpha_1},\ldots,D_t^{1-\alpha_N}\Big),
$
and $D_t^{1-\alpha_i}, i=1,\ldots,N$ are the Riemann-Liouville derivatives; the factors $K_{\alpha_i}=\sigma^2/B_{\alpha_i}$ of the diagonal matrix represent diffusion coefficients with the dimension ${\rm cm}^2/{\rm sec}^{\alpha}$; it can be noted that if $N=1$ the usual fractional Fokker-Planck equation \cite{metzler:00} is recovered.

Next our aim is to calculate the PDF of finding the particle at position $x$ at time $t$, denoted as $g(x,t)$, and the mean squared displacement (MSD) of the process. Without loss of generality, we assume $0<\alpha_1\leqslant \alpha_2\leqslant\ldots\leqslant\alpha_N<1$, and take the equilibrium and initial distributions, respectively, as $\big|\rm eq_M\big>=(\varepsilon_1,\ldots,\varepsilon_N)^{T}$ and $\big<{\rm init}\big|=(\lambda_1,\ldots,\lambda_N)$. We denote the matrix $I-M^{T}\Phi(s)\Lambda(k)$ in Eq. \eqref{eq0} as $A(s)$; it is irreducible or not if and only if $M$ is or not. Plugging the distributions of waiting time and jump length into $A(s)$, for the irreducible transition matrix $M$ the asymptotic expression of the inverse matrix of $A(s)$ can be expressed as
%Now we consider the matrix $I-M^{T}\Phi(s)\Lambda(k)$ in Eq. \eqref{eq0} and denote this matrix as $A(s)$. According to whether the transition matrix M reducible or not, we consider 2 cases. If the transition matrix is irreducible, then $A(s)$ is also irreducible. And we take the distributions of waiting time and jump length which have been shown above into matrix $A(s)$. With the help of eigenvalues, left and right eigenvectors, we can obtain the asymptotic behaviour of the inverse matrix of $A(s)$ that is
$$
A^{-1}(s)\thicksim \frac{\big|{\rm eq_M}\big>\big<\Sigma\big|}{\big<\Sigma\big|\Psi(s)\big|{\rm eq_M}\big>+\sigma^2k^2\big<\Sigma\big|{\rm eq_M}\big>}.
$$
Then we obtain
\begin{equation}\label{eq3}
g(k,s)=\big<\Sigma\big|G(k,s)\big>\thicksim\frac{1}{s}\frac{\big<\Sigma\big|\Psi(s)\big|{\rm eq_M}\big>}{\big<\Sigma\big|\Psi(s)\big|{\rm eq_M}\big>+\sigma^2k^2}.
\end{equation}
The PDF in the Fourier-Laplace space given in Eq. \eqref{eq3} is the same as the PDF of (natural-form) distributed-order diffusion discussed in \cite{sandev:00,chechkin:00} with $p(\alpha)=\sum^{N}_{i=1}\varepsilon_i\delta(\alpha-\alpha_i)$, even though their backgrounds are completely different. From Eq. \eqref{eq3}, the MSD can be got as
$$
\big<x^2(t)\big>\sim\mathcal{L}^{-1}\Bigg\{\frac{2\sigma^2}{s(\varepsilon_1B_{\alpha_1}s^{\alpha_1}+\ldots+\varepsilon_NB_{\alpha_N}s^{\alpha_N})}\Bigg\}.
$$
Therefore, when $t$ is large enough ($s$ tends to $0$), the MSD of the process behaves as $\big<x^2(t)\big>\sim\mathcal{L}^{-1}\Big\{\frac{2K_{\alpha_1}}{\varepsilon_1s^{1+\alpha_1}}\Big\}\sim\frac{2K_{\alpha_1}}{\varepsilon_1\Gamma(1+\alpha_1)}t^{\alpha_1}$; i.e., the MSD of the process with the irreducible transition matrix behaves asymptotically as $t^\alpha$, where $\alpha$ is the smallest one among all exponents of the power-law waiting time distributions of the internal states. When the transition matrix of the Markov chain of the internal states is irreducible, the equilibrium distribution $\big|{\rm eq_M}\big>$ does not depend on the initial distribution, and the initial distribution has no influence on the final PDF or MSD.

%Through CTRW model by choosing the distribution of jump length $\lambda(k)\sim 1-\sigma^2k^2$ and the distribution of waiting time $\phi(s)=[1+B_{\alpha_1}\varepsilon_1s^{\alpha_1}+\ldots+B_{\alpha_N}\varepsilon_Ns^{\alpha_N}\big]^{-1}$ we can also obtain Eq.\eqref{eq1}. Besides if time is sufficient long, we will have $\phi(t)\sim\frac{\alpha_1}{\Gamma(1-\alpha_1)}\frac{B_{\alpha_1}\varepsilon_{\alpha_1}}{t^{\alpha_1+1}}$. This is the relation between this process and the CTRW model.

\begin{figure}[!h]
\centering
\subfigure{\includegraphics[height=3cm,width=4cm]{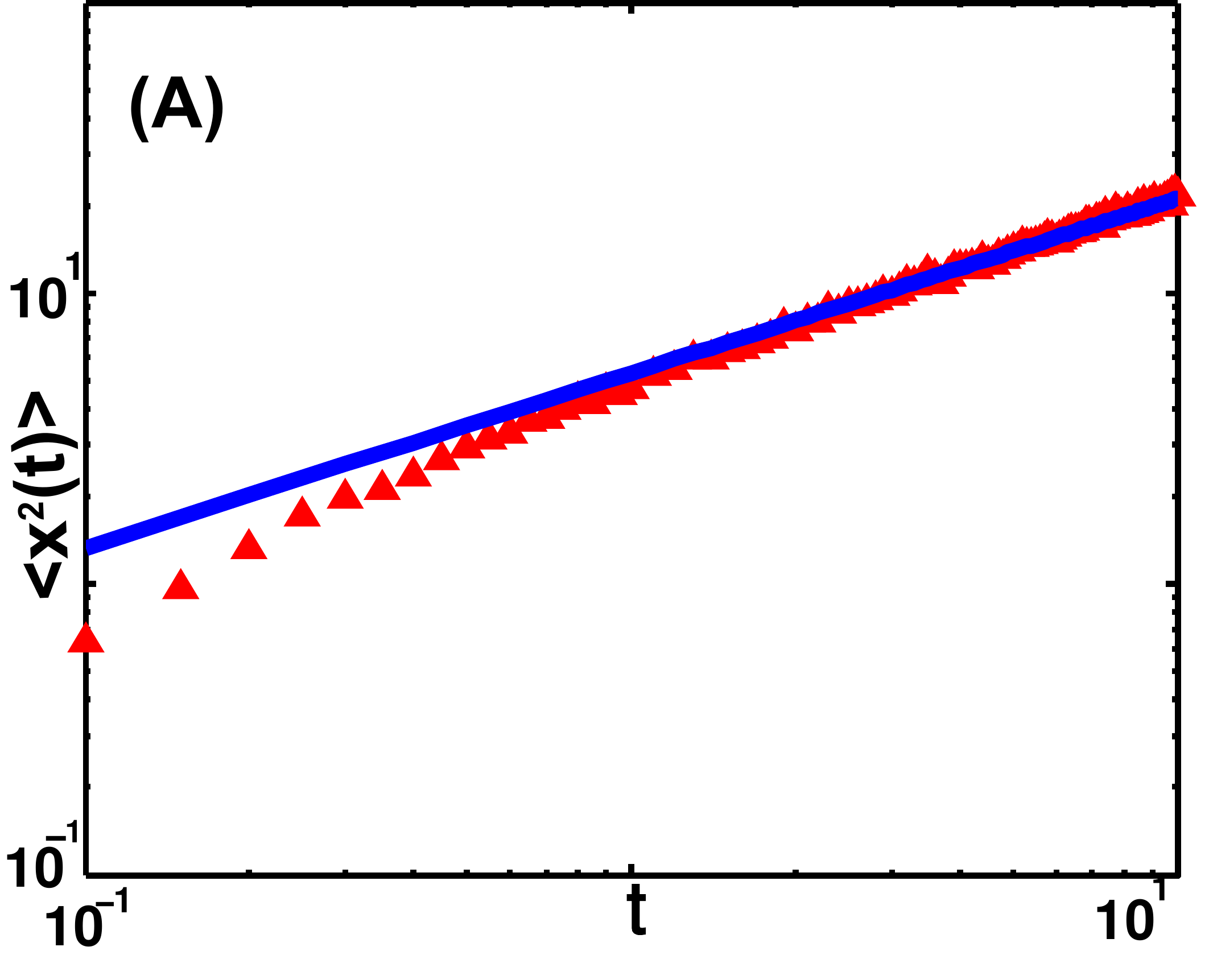}}
\subfigure{\includegraphics[height=3.2cm,width=4cm]{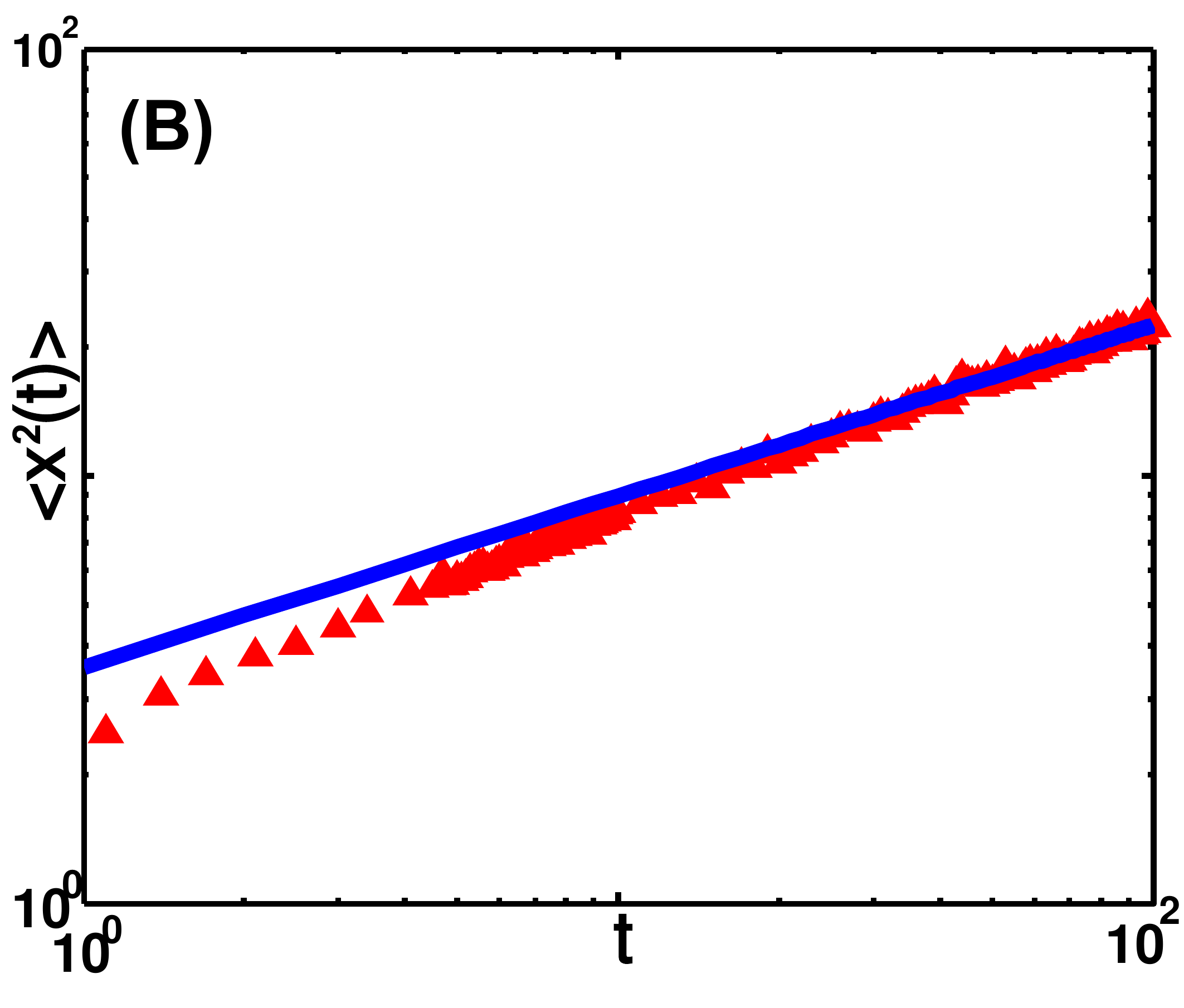}}
\caption{The evolution of the MSD of the process with internal states, sampled over $10^5$ realizations; the solid lines are analytical results, while the triangles are simulation ones. The transition of the internal states forms a reducible Markov chain with the transition matrix $M={\rm diag}\{M_1,M_2,M_3\}$, where $M_1=
\begin{pmatrix}
1/2 & 1/2\\
1/2 & 1/2
\end{pmatrix}$, $M_2=
\begin{pmatrix}
1/3 & 2/3\\
2/3 & 1/3
\end{pmatrix}$, and $M_3=
\begin{pmatrix}
1/4 & 3/4\\
1/2 & 1/2
\end{pmatrix}$.
The exponents $\alpha_1=0.2$, $\alpha_2=0.3$, $\alpha_3=0.4$, $\alpha_4=0.5$, $\alpha_5=0.6$, and $\alpha_6=0.8$. The initial distribution of (A) is $\big|{\rm init}\big>=(1/6,1/6,1/6,1/6,1/6,1/6)^T$, which theoretically implies that the MSD behaves as $t^{0.6}$ (solid line). The initial distribution of (B) is $\big|{\rm init}\big>=(1/4,1/4,1/4,1/4,0,0)^T$, theoretically signifying the evolution of the MSD like $t^{0.4}$ (solid line).
%The solid lines in Fig(A) and Fig.(B) are the theory results. The circles are the simulation results. Both of Fig.(A) and Fig.(B) fit very well. So the initial distribution affects the final results if the process of internal states is reducible.
}%0.6
\label{fig1}
\end{figure}

On the other hand, if the transition matrix is reducible, then the initial distribution of the internal states often makes an important impact on the final results. First, we consider the case that the transition matrix has the form of ${\rm diag}\big\{M_1(s),\ldots M_j(s)\big\}$,
where the matrices $M_i, i=1,\ldots,j$ are irreducible with the dimension of $n_i\times n_i$, and $n_1+n_2+\ldots+n_j=N$. From the form of the transition matrix, one can see that the internal states of the process actually consist of several independent Markov chains with the transition matrices $M_1, M_2,\ldots, M_j$.
%The connection among these Markov chains is the initial distribution, such as the transit matrix is identity matrix, the most special case.
Following the structure of $M$, since $\Psi(s)$ and $I$ are diagonal, we can rewrite $\Psi(s)$ and $I$ as the form consistng of $\Psi_i(s)$ and $I_i$, $i=1,2,\ldots,j$. Then the matrix $A(s)$ has the form of ${\rm diag}\big\{A_1(s),\ldots, A_j(s)\big\}$, where $A_i(s)=I_i-M_i^T\Psi_i(s)+\sigma^2k^2M_i^T$. The vectors can also be rewritten as $\big|{\rm init}\big>=\big(\big|{\rm init}\big>_1,\ldots,\big|{\rm init}\big>_j\big)^T$, $\big|{\rm eq_M}\big>=\big(\big|{\rm eq_M}\big>_1,\ldots,\big|{\rm eq_M}\big>_j\big)^T$, and $\big<\Sigma\big|=\big(\big<\Sigma\big|_1,\ldots,\big<\Sigma\big|_j\big)$. For the convenience of statement, we redefine the subscripts, i.e., let $\Psi_i(s)$, $\big|{\rm init}\big>_i$, and $\big|{\rm eq_M}\big>_i$, etc be consist of $\{B_{\alpha_{ir}}s^{\alpha_{ir}}\}$, $\{\lambda_{i,r}\}$, and $\{\varepsilon_{i,r}\}$ etc, respectively, where $r=1,2,\ldots,n_i$.
%Here we assume every $M_i$ in the matrix $M$ isn't $\{0\}$.
After obtaining the inverse matrix of $A(s)$, we have the PDF $g(x,t)$ in the Fourier-Laplace space
\begin{equation} \label{PDFReducible}
g(k,s)\sim\sum_{\mathop{i=1} \limits_{i \neq i_1,i_2,\cdots,i_j}}^j \frac{1}{s}\frac{\big<\Sigma_i\big|\Psi_i(s)\big|{\rm eq_M}\big>_i\big<\Sigma_i\big|{\rm init}\big>_i}{\big<\Sigma_i\big|\Psi_i(s)\big|{\rm eq_M}\big>_i+\sigma^2k^2\big<\Sigma_i\big|{\rm eq_M}\big>_i},
\end{equation}
where $i_m=m$ if $\big|{\rm init}\big>_m=0$ otherwise $i_m=0$, $m=1,2,\cdots,j$.
In a very particular case, i.e., the transition matrix is an identity matrix and $i_m=0$ for  $m=1,2,\cdots,N$, the PDF of the process has the form
$$
g(k,s)\sim\sum_{i=1}^N\frac{1}{s}\frac{\lambda_is^{\alpha_i}}{s^{\alpha_i}+K_{\alpha_i}k^2}\sim\frac{1}{s}-k^2\Bigg[\sum_{i=1}^NK_{\alpha_i}\lambda_i s^{-\alpha_i-1}\Bigg],
$$
which is the same as the PDF of (modified-form) distributed-order diffusion \cite{sandev:00,chechkin:00} with $p(\alpha)=\sum^{N}_{i=1}\varepsilon_i\delta(\alpha-\alpha_i)$ (their physical backgrounds are completely different).
%Different from the first case, if the transition matrix is reducible, the initial distribution will affect the equilibrium distribution and the final PDF. Thus in this case, the initial distribution play an important role. Now we consider a very special process whose transition matrix is identity matrix. Then we have the PDF of this process yields in the form
%$$
%g(k,s)\sim\sum_{i=1}^j\frac{1}{s}\frac{\lambda_is^{\alpha_i}}{s^{\alpha_i}+K_{\alpha_i}k^2}\sim\frac{1}{s}-k^2\Bigg[\sum_{i=1}^jK_{\alpha_i}\lambda_i s^{-\alpha_i-1}\Bigg]
%$$
%This result is the same with modified-form distributed-order diffusion with the same weight function $p(\alpha)$ illustrated in irreducible matrix when t large enough.
From (\ref{PDFReducible}), there exists
$$
\big<x^2(t)\big>\sim\mathcal{L}^{-1}\Bigg\{\sum_{\mathop{i=1} \limits_{i \neq i_1,i_2,\cdots,i_j}}^j\frac{2\sigma^2\big<\Sigma_i\big|{\rm init}\big>_i\big<\Sigma_i\big|{\rm eq_M}\big>_i}{s\big<\Sigma_i\big|\Psi_i(s)\big|{\rm eq_M}\big>_i}\Bigg\},
$$
which behaves as $\big<x^2(t)\big>\sim t^{\alpha^\ast}$ for large $t$ with  $\alpha^\ast=\max\limits_{ \mathop{1\leqslant i\leqslant j} \limits_{i \neq i_1,i_2,\cdots,i_j}}\Big\{\min\limits_{1\leqslant r\leqslant n_i}\{\alpha_{ir}\}\Big\}$, being confirmed by simulating the stochastic process (see Fig. \ref{fig1}); the influence of the distribution of the initial states is also observed.
 % And the theoretical results and the simulation fit very well. Under this circumstance, when time is sufficient long, this process can consider as a CTRW model with distribution of the waiting time of the form $\phi(t)\sim\frac{1}{t^{\alpha_{n_j1}+1}}$.

Next, we consider the case that the transition matrix is not strictly the form of block diagonal matrix. Without loss of generality, we assume that the first $i$ rows of the transition matrix still keep the form of block diagonal matrix while the others not. If the elements of the initial distribution from $(i+1)$-th to $N$-th are $0$ (denoting this part of the vector as $\big|{\rm init}\big>_{i+1,N}$), then the results illustrated by Fig. \ref{fig1} still hold, just neglecting the last $(N-i)$ internal states; on the other hand, if $\big|{\rm init}\big>_{i+1,N}$ does not equal to zero and, at the same time, the elements of $\big|{\rm init}\big>$ from first to $i$-th are nonzero, the results for large $t$ can still be displayed by Fig. \ref{fig1}, since all the particles in the states ($(i+1),\cdots, N$) finally go to some of the first $i$ states. A little bit complex case is that as least one of the first $i$ elements of \big|{\rm init}\big> is zero and $\big|{\rm init}\big>_{i+1,N} \neq 0$; in this case, taking $\big<{\rm init}\big|M^{N-1}$ as the new initial distribution, the results are still depicted by Fig. \ref{fig1}, ignoring the last $(N-i)$ internal states.

 %we can still use the above conclusions by neglecting the last (n-i)th internal states. Second, if every factor of the initial distribution isn't 0, because of the absorbing of the first ith states, the last (N-i)th have little affect on the final PDF. Thus we can conclude the MSD of this process have the same form with the case of block diagonal matrix by qualitative analysis. And we only compare the powers of the first ith waiting time distributions to obtain $\alpha_0$. Last, if there exists 0 among the first ith factors of the initial distribution and the last (N-i)th factors are not all zero. Then we can calculate the vector $\big<{\rm init}\big|M^{N-1}$ as the new initial distribution. If the first ith factors of this vector still have 0, then we can neglect the corresponding internal states and the last (N-i)th internal states. And we can treat this process as the previous case. Further if we want to analyze the last 2 cases quantitatively, we can use the equilibrium distribution as the equivalent initial distribution. And we can consider this process as the first case when time is long enough.

\emph{Equations governing the distribution of the functionals of the paths and internal states of the process.}---Two types of functionals will be considered. One is still defined as $A=\int_0^t U(x(\tau)) d\tau$ \cite{Kac:49}, being widely discussed \cite{carmi:00,turgeman:00,cairoli:00,satya:00,Wu:16}, where $U(x)$ is a prescribed function and $x(t)$ is a trajectory of a particle. The other functional is first introduced here, defined as
$A_s=\int^t_0 U(i(\tau))d\tau$, where $i(\tau)$ represents that the particle is in the state $i$ at time $\tau$, naturally its values belonging to $\{1,2,\ldots,N\}$.

%We consider the Feynman-Kac equation of the process. The fractional Feynman-Kac equation for functionals of anomalous paths is introduced by Ref.\cite{carmi:00}.

\begin{figure}[!h]
\centering
\subfigure{\includegraphics[height=3cm,width=4cm]{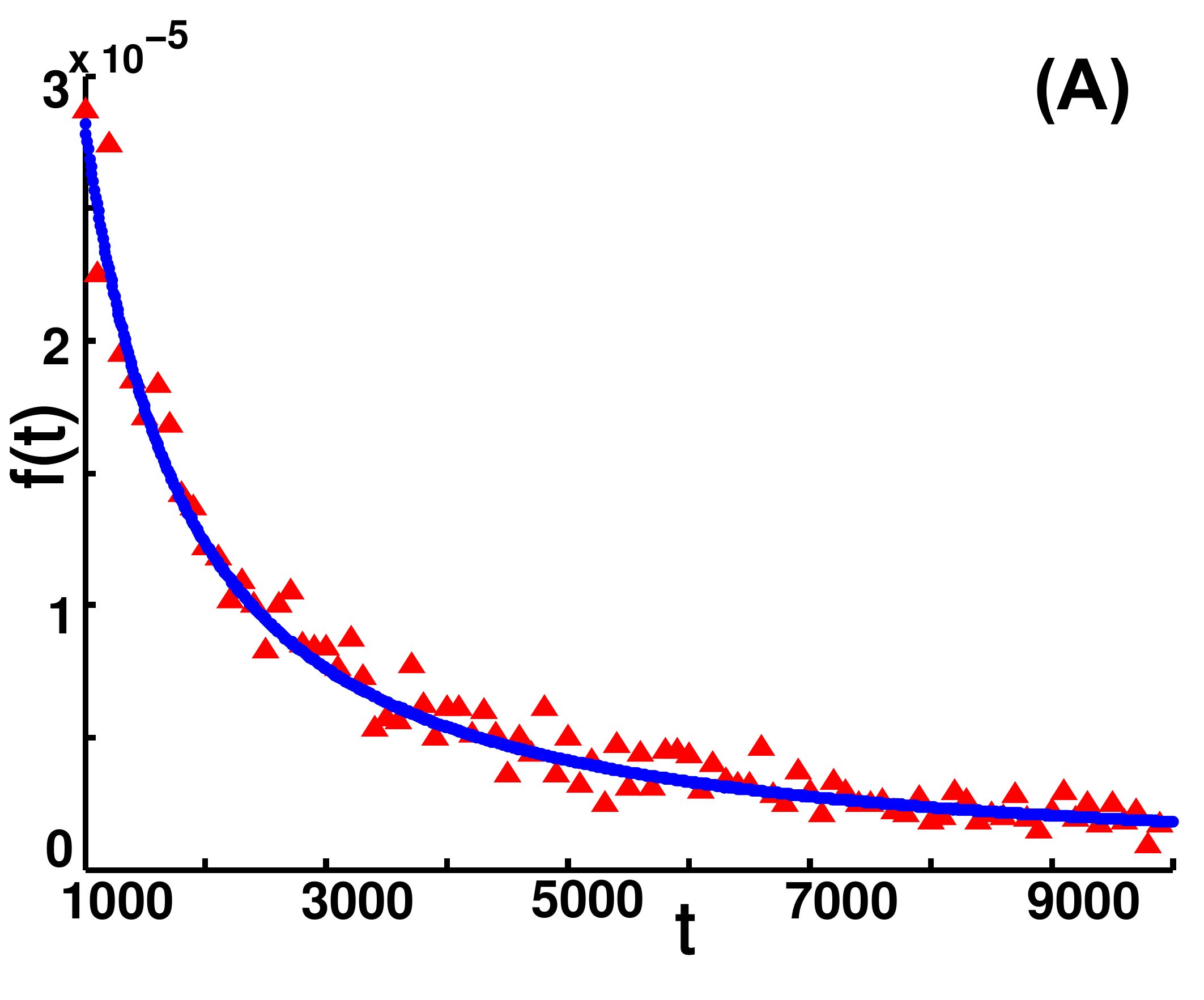}}
\subfigure{\includegraphics[height=3cm,width=4cm]{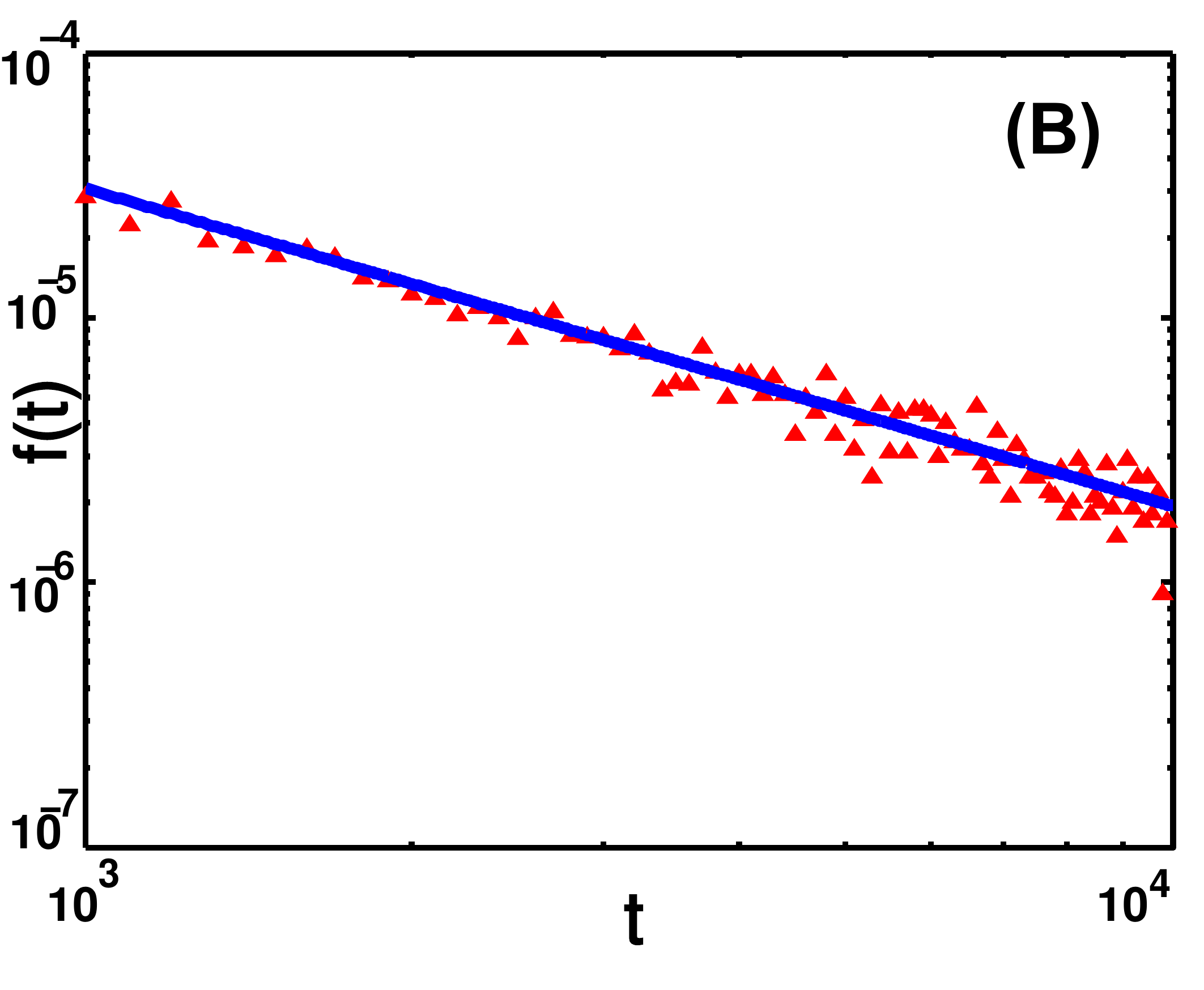}}
\caption{PDF of the first passage time with alternating transition matrix. (A) is for original scale and (B) for the log-log scale with the slope $-1.2$; the real lines are for theoretical results and the triangles for the simulation ones, sampled over $10^5$ realizations. The initial distribution is $\big|{\rm init}\big>=(3/4,1/4)^T$. The other parameters are, respectively, taken as, $\alpha_1=0.8$, $\alpha_2=0.4$, and $B=1$
%And the parameters of the waiting time distributions are chosen to be $\alpha_1=0.4$, $\alpha_2=0.8$. The initial distribution of internal states is $\big|{\rm init}\big>=(1/4,3/4)^T$. We take the parameter B is 1. Fig.(A) is shown in the real time scale, while fig.(B) is illustrated in log-log scale. The triangles in fig.(A) and (B) are the simulation results. The lines are the analysis results. We can see the simulation results and the theory conclusions fit with each other.
}%0.6
\label{fig2}
\end{figure}

We use the notation $g^{(i)}(x,A,t), i=1,2,\ldots,N$ to represent the joint PDF of finding the particle at position $x$ with the functional $A$ and in the internal state of $E_i$ at time t. Gather all $g^{(i)}(x,A,t)$ to form a column vector denote by $\big|G(x,A,t)\big>$. Following the process of the derivation of the fractional Feynman-Kac equation \cite{carmi:00}, we have
\begin{widetext}
%\begin{equation}\label{eq4}
%\big|G(k,\rho,s)\big>=\frac{\Psi(s+\rho U(-i\frac{\partial}{\partial k}))}{s+\rho U(-i\frac{\partial}{\partial k})}\Big[I-M^T+\frac{k^2a^2}{2}M^T+M^T\Psi\Big(s+\rho U\big(-i\frac{\partial}{\partial k}\big)\Big)\Big]^{-1}\big|{\rm init}\big>.
%\end{equation}
%After performing the inverse Fourier and Laplace transform, there exists
\begin{equation}\label{ForwardFeymannKac}
\begin{split}
M^T\frac{\partial}{\partial t}\big|G(x,\rho,t)\big>
=&(M^T-I){\rm diag}\big(B_{\alpha_1}^{-1},\ldots,B_{\alpha_N}^{-1}\big)\mathcal{D}_t^{{\rm diag}(1-\alpha_1,\ldots,1-\alpha_N)}\big|G(x,\rho,t)\big>\\
&+M^T\frac{\partial^2}{\partial x^2}{\rm diag}(K_{\alpha_1},\ldots,K_{\alpha_N})
\mathcal{D}_t^{{\rm diag}(1-\alpha_1,\ldots,1-\alpha_N)}\big|G(x,\rho,t)\big>-\rho U(x)M^T\big|G(x,\rho,t)\big>,
\end{split}
\end{equation}
\end{widetext}
with the initial condition $\big|G(x,A,t)\big>\,\big|_{t=0}=\delta(A)\delta(x)\delta(t)\big|{\rm init}\big>$,  and 
%the operator $\mathcal{D}_t^{{\rm diag}(1-\alpha_1,\ldots,1-\alpha_N)}$ is defined as
$
\mathcal{D}_t^{{\rm diag}(1-\alpha_1,\ldots,1-\alpha_N)}={\rm diag}(\mathcal{D}_t^{1-\alpha_1},\ldots,\mathcal{D}_t^{1-\alpha_N})
$
with $\mathcal{D}_t^{1-\alpha_i}$ being the fractional substantial derivative \cite{carmi:00,friedrich:00}. Next we derive the backward version of Eq. (\ref{ForwardFeymannKac}). We use the notation $\big|{\rm init}\big>_{x_0}=\big(\lambda^{(1)}_{x_0},\lambda^{(2)}_{x_0},\ldots,\lambda^{(N)}_{x_0}\big)$ to represent the initial distribution of the internal states of the process starting at $x_0$, and $g_{x_0}^{(i)}(A,t)$ the PDF of the functional $A$ of the process at $t$, starting at $x_0$ with the internal state $E_i$.
%Thus $\lambda^{(i)}_{x_0}g_{x_0}^{(i)}(A,t)$ represents the same meaning except the probability that the process starts at the $ith$ state. Here we only concentrate on the equations of $g_{x_0}^{(i)}(A,t)$.
After some calculations, one can get
\begin{widetext}
\begin{equation}\label{eq5}
\begin{split}
M\frac{\partial}{\partial t}\big|G_{x_0}(\rho,t)\big>=&{\rm diag}\big(B_{\alpha_1}^{-1},\ldots,B_{\alpha_N}^{-1}\big)\mathcal{D}_t^{{\rm diag}(1-\alpha_1,\ldots,1-\alpha_N)}(M-I)\big|G_{x_0}(\rho,t)\big>\\
&+{\rm diag}(K_{\alpha_1},\ldots,K_{\alpha_N})
\mathcal{D}_t^{{\rm diag}(1-\alpha_1,\ldots,1-\alpha_N)}M\frac{\partial^2}{\partial x_0^2}\big|G_{x_0}(\rho,t)\big>
-\rho U(x_0)M\big|G_{x_0}(\rho,t)\big>.
\end{split}
\end{equation}
\end{widetext}
%We note that the matrices and the operators in Eq.\eqref{eq5} can't exchange.
If one is only interested in the PDF of $A$ at $t$ of the process starting at $x_0$, just calculate $g_{x_0}(A,t)=\sum_{i=1}^N\lambda_{x_0}^{(i)} g^{(i)}_{x_0}(A,t)$. Next, we give a specific application of (\ref{eq5}) for calculating the distribution of the first passage time $t_f$, being the time that the particle starting at $x_0$ ($<B$) first reaches $x=B$.
%When the particles start at $x_0$, we called the time that the particles first hit $x=B$ the first passage time denoted as $t_f$. Specifically,
Define $A_f=\int_0^t U(x(\tau))d\tau$, where $U(x)=0$ if $x<B$ otherwise it equals to $1$.
%$$
%U(x)=
%\begin{cases}
%0 & \text{$x<B$}\\
%1 & \text{$x>B$}
%\end{cases}
%$$
According to \cite{carmi:00,Wu:16}, there exists the relation $Pr\{t_f>t\}=Pr\Big\{\max_{0\leq\tau<t}x(\tau)<B\Big\}=\lim_{\rho\rightarrow\infty}g_{x_0}(\rho,t)$.
For the process of two internal states with the alternating transition matrix, and the coefficients $K_{\alpha_1}=K_{\alpha_2}=B_{\alpha_1}=B_{\alpha_2}=1$,
%Performing Laplace transform on (\ref{eq5}) leads to: for $x_0<B$,
%\begin{widetext}
%\begin{displaymath}
%\begin{cases}
%sg^{(2)}_{x_0}(\rho,s)-1=&-g^{(1)}_{x_0}(\rho,s)s^{1-\alpha_1}+g^{(2)}_{x_0}(\rho,s)s^{1-\alpha_1}+s^{1-\alpha_1}\frac{\partial^2}{\partial x_0^2}g^{(2)}_{x_0}(\rho,s) \\
%sg^{(1)}_{x_0}(\rho,s)-1=&g^{(1)}_{x_0}(\rho,s)s^{1-\alpha_2}-g^{(2)}_{x_0}(\rho,s)s^{1-\alpha_2}+s^{1-\alpha_2}\frac{\partial^2}{\partial x_0^2}g^{(1)}_{x_0}(\rho,s);
%\end{cases}
%\end{displaymath}
%for $x_0>B$,
%\begin{displaymath}
%\begin{cases}
%&sg^{(2)}_{x_0}(\rho,s)-1=-g^{(1)}_{x_0}(\rho,s)(s+\rho)^{1-\alpha_1}+g^{(2)}_{x_0}(\rho,s)(s+\rho)^{1-\alpha_1}+(s+\rho)^{1-\alpha_1}\frac{\partial^2}{\partial x_0^2}g^{(2)}_{x_0}(\rho,s)-\rho g^{(2)}_{x_0}(\rho,s)\\
%&sg^{(1)}_{x_0}(\rho,s)-1=g^{(1)}_{x_0}(\rho,s)(s+\rho)^{1-\alpha_2}-g^{(2)}_{x_0}(\rho,s)(s+\rho)^{1-\alpha_2}+(s+\rho)^{1-\alpha_2}\frac{\partial^2}{\partial x_0^2}g^{(1)}_{x_0}(\rho,s)-\rho g^{(1)}_{x_0}(\rho,s).
%\end{cases}
%\end{displaymath}
%\end{widetext}
we have
%, and demanding $g^{(1)}_{x_0}(\rho,s)$ and $g^{(2)}_{x_0}(\rho,s)$ are finite when $\mid x_0\mid\rightarrow \infty$, besides demanding $g_{x_0}(\rho,s)$ and its first derivative are continue at $x_0=B$. Therefore we have
$$
\lim_{\rho\rightarrow\infty}g_0(\rho,s)=\frac{1}{s}\Bigg[1-{\rm exp}\Big(-\sqrt{\frac{a_0+b_0}{2}}B\Big)\Bigg],
$$
where $a_0=-2+s^{\alpha_1}+s^{\alpha_2}$ and $b_0=\sqrt{4+s^{2\alpha_1}+s^{2\alpha_2}-2s^{\alpha_1+\alpha_2}}$; here $x_0$ is taken as $0$ but not essential. Thus one can obtain the PDF of the first passage time %, denoted as $f(t)$
$$
f(t)=\mathcal{L}^{-1}\Bigg\{\exp\Big(-\sqrt{\frac{a_0+b_0}{2}}B\Big)\Bigg\}.
$$
When $t$ is big and $\alpha_1>\alpha_2$, then $a_0\sim-2+s^{\alpha_2}$ and  $b_0\sim\sqrt{4+s^{2\alpha_2}}\sim2+\frac{1}{4}s^{2\alpha_2}$. Thus, $\exp\Big(-\sqrt{\frac{a_0+b_0}{2}}B\Big)\sim \exp\Big(-\sqrt{\frac{s^{\alpha_2}}{2}}B\Big)\sim 1-\frac{B}{\sqrt{2}}s^{\frac{s^{\alpha_2}}{2}}$, i.e., $f(t)\sim\frac{B}{\sqrt{2}\mid\Gamma(-\alpha_2/2)\mid}t^{-\alpha_2/2-1}$. The result is confirmed by simulations given in Fig. \ref{fig2}.
%Thus we can obtain $a+b\sim s^{\alpha_2}$ i.e. $exp\Big(-\sqrt{\frac{a_0+b_0}{2}}B\Big)\sim exp\Big(-\sqrt{\frac{s^{\alpha_2}}{2}}B\Big)\sim 1-\frac{B}{\sqrt{2}}s^{\frac{s^{\alpha_2}}{2}}$
%And after inverse Laplace transform with respect to s, we have $f(t)\sim\frac{B}{\sqrt{2}\mid\Gamma(-\alpha_2/2)\mid}t^{-\alpha_2/2-1}$. Similarly, when $\alpha_1<\alpha_2$ we have $f(t)\sim\frac{B}{\sqrt{2}\mid\Gamma(-\alpha_1/2)\mid}t^{-\alpha_1/2-1}$.
%If $\alpha_1=\alpha_2=\alpha$, then $f(t)\sim\frac{B}{\mid\Gamma(-\alpha/2)\mid}t^{-\alpha/2-1}$, being consistent with the one given in \cite{carmi:00}. It can be noted that the initial distribution of the internal states does not influence the distribution of the first passage time, since the transition matrix is irreducible.
%
%For $\alpha_1=\alpha_2=\alpha$, we can easily derive $a_0+b_0=2s^\alpha$. Then we have $f(t)\sim\frac{B}{\mid\Gamma(-\alpha/2)\mid}t^{-\alpha/2-1}$. Actually when $\alpha_1=\alpha_2$, the process only has one state, so the result of the last case in accordance with \cite{carmi:00}. From the results, we can see the initial distribution of the internal states has no affect on the distribution of first passage time for this example. This is due to the transit matrix is irreducible, the initial distribution doesn't influent the PDF. Thus this conclusion is inevitable. The simulation of these results are illustrated in Fig.2.

We then turn to the distribution of $A_s=\int^t_0 U(i(\tau))\tau$, ignoring the position $x$ of the particle. Denote $G^{(i)}(A_s,t)$ as the joint PDF of finding the particle with the functional $A_s$ in the state $i$ at time $t$. Let $
\big|G(A_s,t)\big>=\begin{pmatrix}
G^{(1)}(A_s,t),\ldots,G^{(N)}(A_s,t)
\end{pmatrix}^T
$.
Then we get the governing equation
\begin{widetext}
%\begin{equation}\label{eq6}
%\begin{split}
%\big|G(\rho_s,s)\big>=&{\rm diag}\Big\{\frac{[1-\phi_1(s+\rho_sU(1))]}{s+\rho_s\tau U(1)},\ldots,\frac{[1-\phi_1(s+\rho_sU(N))]}{s+\rho_s\tau U(N)}\Big\}\\
%&\times\Big[I-M^T {\rm diag}\{\phi_1(s+\rho_sU(1)),\ldots,\phi_N(s+\rho_sU(N))\} \Big]^{-1}\big|{\rm init}\big>
%\end{split}
%\end{equation}
%That is
\begin{equation}\label{eq6}
M^T\frac{\partial}{\partial t}\big|G(\rho_s,t)\big>=(M^T-I){\rm diag}\Big\{\frac{1}{B_{\alpha_1}}\mathfrak{D}_t^{1-\alpha_1},\ldots,\frac{1}{B_{\alpha_N}}\mathfrak{D}_t^{1-\alpha_N}\Big\}\big|G(\rho_s,t)\big>-\rho_sM^T{\rm diag}\{U(1),\ldots,U(N)\},
\end{equation}
where
\begin{displaymath}
\begin{split}
\mathfrak{D}_t^{1-\alpha_i}G(\rho_s,t)=\frac{1}{\Gamma(\alpha_i)}\Big[\frac{\partial}{\partial t}+\rho_sU(i)\Big]\int_0^t\frac{\exp(-(t-\tau)\rho_sU(i))}{(t-\tau)^{1-\alpha_i}}G(\rho_s,\tau)d\tau.
\end{split}
\end{displaymath}
\end{widetext}
A direct application of (\ref{eq6}) is to calculate the distribution of the fraction of the occupation time, i.e., the distribution of $t^{(i)}/t$, denoted as $l_{t^{(i)}/t}(x)$, where $t^{(i)}$ represents the occupation time of state $i$. Without loss of generality, we only consider the occupation time of the first state by letting $U(i(\tau))=1$ if $i(\tau)=1$, otherwise $U(i(\tau))=0$. Here we just present three results (in the case that the transition matrix is irreducible \cite{reducible:17}): 1) if $\alpha_1<\alpha_2 \leq \alpha_3 \leq \cdots \leq \alpha_N$, then $l_{t^{(i)}/t}(x) \sim \delta(x-1)$; 2) if $\alpha_1$ is not the strictly smallest one, then $l_{t^{(i)}/t}(x) \sim \delta(x)$; 3) if $\alpha_1=\alpha_2=\cdots=\alpha_m$ ($ 2\leq m \leq N$) but smaller than other exponents,  then
\begin{widetext}
$$
\lim_{t\rightarrow\infty}l_{t^{(1)}/t}(x)=\frac{\sin(\pi\alpha)}{\pi}\frac{\varepsilon_1(\varepsilon_2+\ldots+\varepsilon_m)(1-x)^{\alpha-1}x^{\alpha-1}}{(\varepsilon_2+\ldots+\varepsilon_m)^2x^{2\alpha}+\varepsilon_1^2(1-x)^{2\alpha}+2\varepsilon_1(\varepsilon_2+\ldots+\varepsilon_m)\cos(\alpha\pi)x^{\alpha}(1-x)^{\alpha}},
$$
\end{widetext}
from which the classic arcsine law can be recovered \cite{stefani:00}.

\emph{Conclusion.}---We derive the Fokker-Planck equations as well as equations of the functionals of the paths and internal states of the fcP process with multiple internal states. Based on the equations, the MSD is analyzed, and the first passage time and fraction of occupation time are calculated. If the transition matrix of the internal states is reducible, the initial distribution of the internal states significantly influence the final results.

%We have obtained the Fokker-Planck equations as well as equations of the functionals of the paths and internal states of the FCP process with multiply internal states. From these equations we can calculate the PDF and MSD of the FCP processes. We build a bridge between this process with the distributed-order diffusion, and this process with the CTRW model. We also see the initial distribution can only influent the process with reducible transit matrix. And we also obtain the Feynamn-Kac equation, by choosing different function, we can calculate many different values. In this letter, we calculate the distribution of the first passage time. In the final part, we neglect the variable x, only consider fractional Poisson process, we can also obtain its equations. And from this equations, we calculate the distribution of the friction of occupation time to the total time. The result is determined by the transit matrix, the equilibrium distribution, and the relationship among $\alpha_1,\ldots,\alpha_N$. Besides if the transit matrix is reducible, then the initial distribution will also influent the results.


\begin{thebibliography}{10}% Produces the bibliography via BibTeX.

\bibitem{Kleinrock:76} L. Kleinrock, \emph{Queueing Systems: Theory} (John Wiley \& Sons, Canada,  1976).

\bibitem{godreche:00} C. Godr\`{e}che and J. M. Luck, J. Stat. Phys. {\bf 104}, 489 (2001).

\bibitem {laskin:00} N. Laskin, Commun. Nonlinear Sci. Numer. Simul. {\bf 8}, 201 (2003).

\bibitem{burov:00} S. Burov and E. Barkai, Phys. Rev. Lett. {\bf 107}, 170601 (2011).
%laskin: fractal Poisson porcess.
\bibitem {lowen:00} S. B. Lowen and M. C. Teich, Phys. Rev. E {\bf 47}, 922 (1993).


\bibitem {Godec:17} A. Godec and R. Metzler, J. Phys. A: Math. Theor. {\bf 50}, 084001 (2017).


\bibitem{niemann:00} M. Niemann, E. Barkai, and H. Kantz, Math. Model. Nat. Phenom. {\bf 11}, 191  (2016).

\bibitem{meerschaert:00} M. M. Meerschaert and A. Sikorskii, \emph{Stochastic Models for Fractional Calculus} (Walter de Gruyter, Berlin, 2012).

\bibitem {cartea:00} ${\rm\acute{A}}$. Cartea and D. del-Castillo-Negrete, Phys. Rev. E {\bf 76}, 041105 (2007).

\bibitem {metzler:00} R. Metzler and J. Klafter, Phys. Rep. {\bf 339}, 1 (2000).

\bibitem {Metzler:16} R. Metzler, J. -H. Jeon, A. G. Cherstvy and E. Barkai, Phys. Chem. Chem. Phys., {\bf 16}, 24128 (2014).

\bibitem{klafter:00} J. Klafter and I. M. Sokolov, \emph{First Steps in Random Walks: From Tools to Applications} (Oxford University Press, Oxford, 2011).

%lowen: fractal renewal processes generate 1/f noise
\bibitem {golding:00} I. Golding and E. C. Cox, Phys. Rev. Lett. {\bf 96}, 098102 (2006).



%\bibitem {stanislavsky:00} A. Stanislavsky, K. Weron, and A. Weron, Phys. Rev. E {\bf 78}, 051106 (2008).

\bibitem{scalas:00} E. Scalas, Lecture Notes in Econom. and Math. Systems {\bf 567}, 3 (2006).



\bibitem{carmi:00} S. Carmi, L. Turgeman, and E. Barkai, J. Stat. Phys. {\bf 141}, 1071  (2010).

\bibitem{turgeman:00} L. Turgeman, S. Carmi, and E. Barkai, Phys. Rev. Lett. {\bf 103}, 190201 (2009).

\bibitem{cairoli:00} A. Cairoli and A. Baule, Phys. Rev. Lett. {\bf 115}, 110601 (2015).

\bibitem{satya:00} S. N. Majumdar, Curr. Sci. {\bf 89}, 2076 (2005).

\bibitem{feller:00} W. Feller, \emph{An Introduction to Probability Theory and Its Application} (Vol. 1, John Wiley \& Sons, US, 1968).

\bibitem{redner:00} S. Redner, \emph{ A Guide to First-Passage Processes} (Cambridge University Press, Cambridge, 2001).

\bibitem{Deng:17} W. H. Deng, X. C. Wu, W. L. Wang, EPL {\bf 117}, 10009 (2017).

\bibitem {Walker:11} S. G. Walker, Linear Multilinear Algebra {\bf 59}, 755 (2011).


%Stephen G. Walker
%Bounds for the second largest eigenvalue of a transition matrix
%2011 Linear Multilinear Algebra, 59:7, 755-760


\bibitem{sandev:00} T. Sandev, A. V. Chechkin, N. Korabel, H. Kantz, I. M. Sokolov, and R. Metzler, Phys. Rev. E {\bf 92}, 042117 (2015).

\bibitem{chechkin:00} A. V. Chechkin, R. Gorenflo, and I. M. Sokolov, Phys. Rev. E {\bf 66}, 046129 (2002).

\bibitem{Kac:49} M. Kac, Trans. Amer. Math. Soc. {\bf 65}, 1 (1949).

\bibitem{friedrich:00} R. Friedrich, F. Jenko, A. Baule, and S. Eule, Phys. Rev. Lett. {\bf 96}, 230601 (2006).

\bibitem{Wu:16} X. C. Wu, W. H. Deng, and E. Barkai, Phys. Rev. E {\bf 93}, 032151 (2016).

%    XC Wu, WH Deng , E Barkai, Tempered fractional Feynman-Kac equation: Theory and examples, Physical Review E, 93(3), 032151, 2016

\bibitem{stefani:00} F. D. Stefani, J. P. Hoogenboom, and E. Barkai, Phys. Today {\bf 62}, 34 (2009).

\bibitem{reducible:17} For the reducible case, the results will be presented in the other publication.



\end{thebibliography}
\end{document}